\documentclass[10pt]{amsart}
\usepackage{latexsym}
\usepackage{amssymb}
\usepackage{amscd}
\usepackage{amsmath}

\newtheorem{proposition}{Proposition}[section]
\newcommand{\Z}{\mathbb{Z}}
   
\begin{document}

\title{A cusp singularity with no Galois cover by 
a complete intersection}
\author{David E. Anderson}
\thanks{Supported by the NSF's VIGRE Fellowship through the Columbia 
University Department of Mathematics.  The author is greatly indebted 
to Professor Walter Neumann for his help and guidance.}

\begin{abstract}
    With an explicit example, we confirm a conjecture by Neumann and 
    Wahl that there exist cusps with no Galois cover by a complete 
    intersection.  Some computational techniques are reviewed, and a 
    method for deciding whether a given cusp has a 
    complete intersection Galois cover is developed.
\end{abstract}

\maketitle

\section{Introduction}

In \cite{wahl}, Neumann and Wahl prove that the universal abelian 
cover of every quotient-cusp is a complete intersection, and 
conjecture generally that a similar result holds for any 
$\mathbb{Q}$-Gorenstein normal surface singularity whose link is a 
rational homology sphere.  (For discussions of these singularities and 
others known to have complete intersection abelian covers, see 
\cite{n2,nw2,wahl}.  By a \textit{cover} of a complex surface singularity we 
mean a finite cover of a germ of the singularity, branched only at the 
singular point.)  In the same article (Proposition 4.1), it is proved 
that every cusp has a cover by a complete intersection, but there exist 
cusps with no abelian cover by a complete intersection.  The authors 
go on to conjecture that some cusps do not even have a Galois cover by 
a complete intersection.

In the present paper, we exhibit a cusp confirming this latter conjecture:  
the cusp classified by the cycle $(8,2,4,3,12)$ satisfies as an example.  
The existence of such a singularity precludes a natural 
generalization of the main theorem of \cite{wahl} from the case of 
quotient-cusps to a statement about cusps.  More generally,  
Conjecture 1.1 of \cite{wahl} does not naturally generalize to the 
case of singularities whose links are not rational homology spheres.

We begin with a brief review of singularity links and cusps; for more 
details, see \cite{neumann,wahl}.

Consider a germ of a normal complex suface singularity $(V,p)$.  A 
\textit{good resolution of the singularity} is a map $\pi:\overline{V} 
\to V$ such that $\overline{V}$ is non-singular, $\pi
:\pi^{-1}(V - \{p\}) \to V - \{p\}$ is an isomorphism, and
$\pi^{-1}(p)$ is a complex curve whose only singularities are normal
crossings and all of whose components (the \textit{exceptional 
curves}) are nonsingular.  The good resolution with fewest
exceptional curves is unique and is called the \textit{minimal good
resolution}.  The \textit{resolution graph} is the weighted graph whose
vertices stand for the exceptional curves weighted according to their
self-intersection numbers, and whose edges correspond to intersections
of exceptional curves.  A \textit{link of the singularity $(V,p)$} is
the oriented 3-manifold $M$ that forms the boundary of a closed
regular neighborhood of $p$.  Using the plumbing procedure described
in \cite{neumann}, the singularity link $M$ can be reconstructed from
the resolution graph.  The topology of the minimal good resolution
$\overline{V}$ is determined by the oriented manifold $M$ (see
\cite{neumann}).

The \textit{cusp singularities} are those whose resolution graphs are cyclic 
and all of whose exceptional curves are rational.  Equivalently, the link 
of a cusp is a torus bundle over the circle, with 
monodromy $A$ having trace $\geq 3$.  A cusp singularity also has a 
\textit{dual cusp} whose link is $-M$ ($M$ with reversed 
orientation), a torus bundle with monodromy $A^{-1}$.

The plumbing calculus of \cite{neumann} reconstructs the monodromy $A$ 
of a cusp's link as follows:  Up to conjugation, the classifying 
matrix is
\[ A = \left(\begin{array}{cc} -e_{k}&1 \\ -1&0 \end{array}\right) \ldots 
\left(\begin{array}{cc} -e_{1}&1 \\ -1&0 \end{array}\right), \]
where $e_{1},\ldots,e_{k}$ are the weights of the vertices of the 
resolution cycle.  We shall call $(-e_{1},\ldots,-e_{k})$ the 
\textit{resolution cycle} of the link determined by $A$, or more 
concisely, the \textit{cycle} of $A$.

Finally, a germ of a complex surface is a \textit{complete intersection} 
(abbreviated ``CI'') if 
there exists an imbedding $(V,p) \hookrightarrow (\mathbb{C}^{n},0)$ 
such that $V$ is given by $(n-2)$ equations.  
It follows from Karras \cite{karras} that a cusp is a CI if and only if the 
resolution cycle of the dual cusp has length $\leq 4$ (see 
\cite{wahl}).

A useful tool in examining Galois covers of a cusp is the fact, 
Proposition \ref{prop:dual}, that 
if a cusp is a Galois cover of a given cusp, then the cusp dual to the 
covering cusp is also a Galois cover.  Thus a cusp has a Galois cover 
by a CI if and only if it has a Galois cover by a cusp with cycle of 
length at most 4.

The topology of the minimal good resolution of a cusp singularity is determined 
(up to orientation reversal, i.e., replacing the cusp by its dual)
by the fundamental group of the singularity link $M$ 
\cite{hirzebruch,neumann}, so the topological question of whether a 
CI Galois cover exists reduces to an algebraic problem.  The 
fundamental group $\pi_{1}(M)$ is isomorphic to the semidirect 
product $\Z^{2} \rtimes_{A} \Z$; 
Galois coverings correspond to normal subgroups of $\pi_{1}(M)$.  
The question now becomes
``Is there a normal subgroup of $\Z^{2} \rtimes_{A} \Z$ whose 
classifying matrix is conjugate to  $B = \left(\begin{smallmatrix} -e_{k}&1 \\ -1&0 
\end{smallmatrix}\right) \ldots \left(\begin{smallmatrix} -e_{1}&1 \\ -1&0 
\end{smallmatrix}\right)$, for some $e_1,...,e_k$ with $k \leq 4$ and 
each $e_{i} \leq -1$?''

To find and verify an example, we proceed as follows.  In Section 
\ref{sec:alg}, we state and prove some algebraic facts about normal 
subgroups of $\pi_{1}(M)$, and arrive at a condition on the 
factorization of the trace of the monodromy that reduces the number of 
candidates for normal subgroups.  Then, in Section \ref{sec:cycles}, 
we discuss the algorithms developed in \cite{neumann} for calculating 
the plumbing cycle of a cusp given the monodromy of its link, and for 
determining the cycle of the dual cusp given the cycle of the 
original.  Guided by conditions on the trace, we then posit a matrix 
to serve as a candidate for the monodromy of the link of a cusp with 
no Galois cover by a CI.  Finally, in Section 
\ref{sec:classify}, we systematically examine the normal subgroups of the 
fundamental group of this link to verify that they have cycles longer 
than 4.

\section{Algebraic Preliminaries} \label{sec:alg}

The fundamental group of the torus bundle $M$ with monodromy $A$ is 
$\Z^{2} \rtimes_{A} \Z$, where the group operation in the semidirect product 
is given by $(x,n)(y,m) := (x+A^{n}y,n+m)$.  In the following, we 
shall denote the monodromy matrix of $M$ by $A$, and the group 
$\pi_{1}(M)$ by $G$.

In discussing subgroups of $\pi_{1}(M)$, we shall follow the terminology 
of \cite{wahl}: a proper subgroup of
$\Z^{2} \rtimes \Z$ that surjects onto $\Z$ will be called a
\textit{covering in the fiber} (or a \textit{fiberwise cover}), and a 
subgroup that contains $\Z^{2}$ will be a \textit{covering in the base}.

\begin{proposition}  \label{prop:comp}
    Every Galois covering is \\
    (a) a Galois covering in the base followed by a 
    Galois covering in the fiber, and \\
    (b) up to isomorphism, a covering in the fiber (not necessarily Galois) 
    followed by a covering in the base. 
\end{proposition}

\begin{proof}
    (a) Let $G$ be the fundamental group of $M$, and let $N$ be the 
    normal subgroup corresponding to the Galois covering (i.e., the 
    fundamental group of the cover).  Denote the projection 
    $\Z^{2}\rtimes_{A}\Z \to \Z$ by $\pi$, let $L = N \cap \Z^{2}$, 
    and let $n=[\Z:\pi(N)]$.  We can construct the 
    following diagram of short exact sequences:
    
    \begin{eqnarray*} \begin{CD}
        1 @>>>     L         @>>>    N   @>>> n\Z @>>> 1     \\
        &      &   @VVV            @VVV       @VVV    \\     
        1 @>>>   \Z^{2}      @>>>    G   @>>>  \Z @>>> 1  .\\
    \end{CD} \end{eqnarray*}
    
    \noindent We want an $\tilde{N}$ to satisfy 
    $N \triangleleft \tilde{N} \triangleleft G$, with a fiber factor equal 
    to that of $G$ and base factor equal to that of $N$:
    \begin{eqnarray*} \begin{CD}
	1 @>>>     L           @>>>   N       @>>> n\mathbb{Z} @>>> 1     \\
	&   & \bigtriangleup &   &\bigtriangleup &  &   \|  &   \\     
	1 @>>>  \mathbb{Z}^{2} @>>> \tilde{N} @>>> n\mathbb{Z} @>>> 1     \\
	&   &     \|         &   &\bigtriangleup & & \bigtriangleup  &    \\     
	1 @>>>  \mathbb{Z}^{2} @>>>    G      @>>>  \mathbb{Z} @>>> 1 \, .
    \end{CD} \end{eqnarray*}
    All we need to show is $\tilde{N} \triangleleft G$, as $N 
    \triangleleft \tilde{N}$ follows.
    This is easy, though, since $\tilde{N}$ can be taken to be the kernel of 
    the natural homomorphism from $G$ to $\Z/n$ given by 
    $(x,y) \mapsto y \pmod n$.

    (b) Now we want $\tilde{N'}$ to satisfy
    
    \begin{eqnarray*} \begin{CD}
	1 @>>>   L    @>>>   N        @>>> n\Z @>>> 1  \\
	&    &  \|    &    &\cap    &  &  \cap  &   \\     
	1 @>>>   L    @>>> \tilde{N'} @>>> \Z @>>> 1  \\
	&   &   \cap  &    &\cap   &  &    \|  &    \\     
	1 @>>> \Z^{2} @>>>    G       @>>> \Z @>>> 1 \, .
    \end{CD} \end{eqnarray*}

    \noindent $N$ is generated by $L$ and $(x,n)$ for some $ x \in 
    \Z^{2}$.  Any subgroup $\tilde{N'}$ that fits in the above 
    horizontal exact sequence will be generated by $L$ and $(x',1)$ 
    for some $x' \in \Z^{2}$.  This $\tilde{N'}$ will contain $N$ only if 
    $(x,n)$ is in the group generated by $L$ and $(x',1)^{n}$, and an 
    $x'$ with this property cannot always be found.  
    However, $N$ is isomorphic to the subgroup of $\Z^{2}\rtimes_{A}\Z$
    generated by $L$ and $(0,n)$, and if we replace $N$ by this
    subgroup then we may take $\tilde{N'}$ as the subgroup generated by
    $L$ and $(0,1)$.
\end{proof}

Using Proposition \ref{prop:comp}(b), we can restrict our search to 
covers of degree at most four in the base: If $\tilde{N'}$ has 
monodromy $B$, then $N = L \rtimes_{B} n\Z = L 
\rtimes_{B^{n}} \Z$ has monodromy $B^{n}$, which has cycle of length $n$ 
times that of $B$.  Covers of degree greater than four thus have 
cycles of length at least 5.  (This is easy to see, since a covering of 
degree $n$ in the base just makes the cycle repeat $n$ times.)

Proposition \ref{prop:comp}(a) tells us we need only look at 
\textit{normal} covers in the fiber; this condition on $N$, together 
with the following, will also 
limit our search.

\begin{proposition} \label{prop:btwn}
    Let $N = L \rtimes n\Z$ be a subgroup of $G$.  Then N is normal 
    if and only if $AL = L$ and $(A^{n} - I)\Z^{2} \subset L \subset \Z^{2}$.
\end{proposition}

\begin{proof}
    $N$ is generated by $L$ and $(x,n)$ for some $x$, and $L$ is a 
    characteristic subgroup of $N$, so $N$ is normal only if $L$ is 
    normal.  $L \to (L,0)$, so conjugation by $(0,1)$ gives
    \[ (0,1)^{-1} (L,0) (0,1) = (A^{-1}L,0). \]
    $L$ is therefore preserved by conjugation if and only if
    \begin{eqnarray*}
	AL = L.
    \end{eqnarray*}
    
    Next, conjugation of $(x,n)$ by $(b,0)$ for any $b \in 
    \Z^{2}$ gives
    \[ (b,0)^{-1} (x,n) (b,0) = (-b + x + A^{n}b,n), \]
    so we must have
    \begin{eqnarray*}
	(A^{n}-I)\Z^{2} \subset L.
    \end{eqnarray*}
    
    The converse is easily verified.
\end{proof}

The above propositions let us give the proof of the following:

\begin{proposition}\label{prop:dual} 
    If a cusp is a Galois cover of a given cusp then its dual cusp
    also is.
\end{proposition}

\begin{proof}
    Given a cusp with monodromy $A$, a covering in the fiber with fiber
    group $L=(A-I)\Z^{2}$ is called a \textit{discriminant cover} (see
    \cite{wahl}). Let us call a cover of $M$ a \textit{sub-discriminant
    cover} if it lies between $M$ and a discriminant cover (i.e., it is
    covered by a discriminant cover). Then Propositions \ref{prop:comp}
    and \ref{prop:btwn} imply that any Galois cover is the result of a
    cover in the base followed by a sub-discriminant cover in the fiber,
    and that any such cover is a Galois cover. In Section 4 of \cite{wahl} it is
    shown that the collection of sub-discriminant covers of a given cusp
    is closed under taking duals, so the proposition follows.
\end{proof}

Proposition \ref{prop:btwn} also tells us we can restrict the number of 
subgroups of $\Z^{2}$ by keeping $\Z^{2}/(A^{n} - I)\Z^{2}$ 
as simple as possible.  It is easy to see that

\begin{proposition} \label{prop:index}
    $|\Z^{2}/(A^{n} - I)\Z^{2}| = 
    |2 - P_{n}(\mathrm{tr}A)|$, where $P_{n}(\mathrm{tr}A)$ is a polynomial in 
    the trace of $A$.
\end{proposition}

\noindent In particular, we have for $n=1,2,3,4$
\begin{eqnarray}
    \left[\Z^{2}:(A^{1} - I)\Z^{2}\right] &= &x - 2  \nonumber \\
    \left[\Z^{2}:(A^{2} - I)\Z^{2}\right] &= &(x - 2)(x + 2) \nonumber \\
    \left[\Z^{2}:(A^{3} - I)\Z^{2}\right] &= &(x - 2)(x + 1)^{2} \nonumber \\
    \left[\Z^{2}:(A^{4} - I)\Z^{2}\right] &= &x^{2}(x - 2)(x + 2), \label{eqn:index}
\end{eqnarray}
where $x$ is the trace of $A$.  By Proposition \ref{prop:btwn}, these 
are the maximal indices (in the fiber) of normal subgroups.  To 
minimize the number of potentially normal subgroups, 
we minimize the number of factors of $|\Z^{2}/(A^{n} - I)\Z^{2}|$; 
let us require that $x$ and $x-2$ be 
prime, $x+2$ be three times a prime, and $x+1$ be twice a prime.  For convenience, 
let us denote the prime factors of $|\Z^{2}/(A^{n} - I)\Z^{2}|$ as 
follows:
\begin{eqnarray}
    p &= &\mathrm{tr}A \nonumber \\
    q &= &p-2 \nonumber \\
    3r &= &p+2 \nonumber \\
    2s &= &p+1, \label{primes}
\end{eqnarray}
where $p$, $q$, $r$, and $s$ are all prime.

\section{Cycles} \label{sec:cycles}

We now turn our attention to the cycle of a cover.  Recall that the 
resolution graph for the link of a cusp singularity is cyclic, and 
each vertex with weight $e_{i}$ corresponds to a joining of trivial torus 
bundles by matrices of the form 
$\left(\begin{smallmatrix} -e_{i}&1 \\ -1&0 \end{smallmatrix}\right)$.  
The monodromy of the entire link, then, is
\begin{eqnarray}
    B = \left(\begin{array}{cc} b_{k}&1 \\ -1&0 \end{array}\right) \ldots
        \left(\begin{array}{cc} b_{1}&1 \\ -1&0 \end{array}\right), 
        \label{eqn:cycle}
\end{eqnarray}
where $b_{i} = -e_{i}$.  Written this way, $B$ classifies the link of a complete 
intersection (up to orientation) iff $k \leq 4$.  Now we need a method 
for expanding a given matrix in this form and an algorithm for 
finding the corresponding dual matrix.  The tools for both are 
described in \cite{neumann}; we summarize briefly here.

Given $B = \left(\begin{smallmatrix} a&b \\ c&d 
\end{smallmatrix}\right)$ with trace $t$, let
\[ \omega = \frac{a-d + \sqrt{t^{2}-4}}{2b}. \]
To find the cycle, first compute the ``continued fraction with minus 
signs'' expansion of $\omega$:
\[\omega = a_{0}-\frac{1}{a_{1}-\frac{1}{a_{2}-\frac{1}{\ddots}}} \]
and let ${[c_{1},\ldots,c_{l}]}$ be the shortest period of this ultimately 
periodic continued fraction.  Put $C = 
\left(\begin{smallmatrix} c_{l}&1 \\ -1&0 \end{smallmatrix}\right) \ldots 
\left(\begin{smallmatrix} c_{1}&1 \\ -1&0 \end{smallmatrix}\right)$, 
with trace $t'$, and solve
\begin{eqnarray}
    \frac{t + \sqrt{t^{2}-4}}{2} = \left(\frac{t' + 
    \sqrt{(t')^{2}-4}}{2}\right)^{n}
\end{eqnarray}
for $n$.  Then the cycle of $B$ is given by concatenating 
$(c_{1},\ldots,c_{l})$ with itself $n$ times.

Since this method depends on finding the period of a continued 
fraction, computation will be greatly simplified if the continued 
fraction is purely periodic.  The following proposition gives a 
condition for pure periodicity of a continued fraction with minus signs.  
The proof is analogous to that of a similar statement about ordinary 
continued fractions; see for example Theorem 11.5.2 in \cite{kumanduri}.

\begin{proposition} \label{prop:pure}
    The continued fraction (with minus signs) expansion of a 
    quadratic irrationality $x = (P+\sqrt{n})/Q$ is 
    purely periodic if and only if $x>1$ and $0<\overline{x}<1$,
    where $\overline{x}$ is the conjugate of $x$ in $\mathbb{Q}{[\sqrt{n}]}$.
\end{proposition}

In relation to the cycle of $B$, Proposition \ref{prop:pure} requires
\begin{eqnarray}
    \frac{(a-d)+\sqrt{t^{2}-4}}{2b} > 1 \;\;\mathrm{and}\;\;
    \frac{2b}{(a-d)-\sqrt{t^{2}-4}} > 1,
\end{eqnarray}
where $t=a+d$.  For $t \geq 3$ and $a,b,d \in \Z$, this reduces to the condition
\begin{eqnarray}
    a > b > -d \geq 0.\label{pure-entries}
\end{eqnarray}

Since the dual to a Galois cover of a given cusp is also a Galois cover, 
we must inspect the dual cycle as well.
If $(b_{1},\ldots,b_{k})$ is the cycle of $B$, then from 
\cite{hirzebruch} we know $b_{i} \geq 2$ for each $i$, 
and $b_{i} \geq 3$ for some $i$.  Write the cycle as
\begin{eqnarray*}
    (b_{1},\ldots,b_{k}) = 
    (m_{1}+3,\underset{n_{1}}{\underbrace{2,\ldots,2}},m_{2}+3,\ldots,
    m_{s}+3,\underset{n_{s}}{\underbrace{2,\ldots,2}});
\end{eqnarray*}
the dual cycle is then given by
\begin{eqnarray*}
    (d_{1},\ldots,d_{l}) = 
    (n_{s}+3,\underset{m_{s}}{\underbrace{2,\ldots,2}},n_{s-1}+3,\ldots,
    n_{1}+3,\underset{m_{1}}{\underbrace{2,\ldots,2}}).
\end{eqnarray*}
As is well-known, the length of the dual cycle can be 
written simply in terms of the entries of the original cycle as 
$\sum_{i=1}^{n}(b_{i}-2)$. 

For our example, we want to find a matrix whose trace satisfies the primality 
condition (\ref{primes}), and whose entries satisfy the condition for 
immediate periodicity given in (\ref{pure-entries}).

First we choose a trace.  The first four numbers 
satisfying (\ref{primes}) are 5, 13, 
1621, and 6661.  Computation time increases rapidly with the size 
of the trace, but 13 is too small to yield a long enough cycle, so 
we use 1621.

Now we construct a matrix in $SL_{2}(\Z)$ with trace 1621 whose entries 
satisfy $a > b > -d \geq 0$.  Some experimentation suggests
\begin{eqnarray}
    A = \left(\begin{array}{cc} 1640 & 221 \\ -141 & -19 \end{array}\right)
\end{eqnarray}
as an example of a matrix that fails the test for a monodromy 
of a CI: the cycle of $A$ is 
$(8,2,4,3,12)$, and the cycle of its dual has length 19.

\section{Classifying Galois covers} \label{sec:classify}

We now verify that the above choice for $A$ does indeed classify a 
cusp with no CI Galois cover.  As mentioned in Section \ref{sec:alg}, 
we only need to consider covers of degree 1, 2, 3, and 4 in the base.  
In the notation of (\ref{primes}), the maximal indices in $\Z^{2}$ for 
the fiber-subgroups of normal subgroups are $q$, $3rq$, $q(2s)^{2}$, 
and $3p^{2}qr$, respectively.  With our choice of $A$ we have 
$p=1621$, $q=1619$, $r=541$, and $s=811$.

For each degree, we find a basis representation for the subgroups 
lying between $\Z^{2}$ and $(A^{n}-I)\Z^{2}$; by Proposition 
\ref{prop:btwn} the normal subgroups are those preserved by $A$.  Once we 
have found such a subgroup, we need to know how a generator $g$ of $\Z$ 
acts on it; the matrix of this action (raised to the $n$th power, for 
degree $n$) will be the monodromy of the Galois cover.  We 
can then use the machinery described in Section \ref{sec:cycles} to 
compute the classifying cycle of the cover; for computation, the 
algorithms were implemented in {\it Mathematica}.  The 
action of $g$ on the maximal-index subgroup $(A^{n}-I)\Z^{2}$ is always 
by $A$, since $A$ commutes with $(A^{n}-I)$, so we only need to check
subgroups properly intermediate between $\Z^{2}$ and $(A^{n}-I)\Z^{2}$.

In fact, we need to treat only half of the possible subgroups 
for each degree.  The other half represent dual cusps, and we know 
how to compute their cycles from the original cusps.  In general, if 
${[\Z^{2}:(A^{n}-I)\Z^{2}]} = xy$, and $L$ is a subgroup of index $x$, 
then a subgroup $L'$ of index $y$ represents the fiber of the dual cusp.  
(This is implied in Section 4 of \cite{wahl}.)
% (This follows from Proposition 4.2 of \cite{wahl}: Let $L$ 
% represent the fiber of a cusp, and write $K = L/(A^{n}-I)\Z^{2}$.  Then 
% the orthogonal complement of $K$ in the \textit{discriminant group} $D \cong 
% \Z^{2}/(A^{n}-I)\Z^{2}$ can be written $K^{\perp} = 
% L'/(A^{n}-I)\Z^{2}$, where $L'$ represents the dual 
% cusp.  Now ${[\Z^{2}:L]} = |K^{\perp}|$ and 
% ${[\Z^{2}:L']} = |K|$, so the product of these indices is 
% equal to the index ${[\Z^{2}:(A^{n}-I)\Z^{2}]}$.)

\subsection{One-fold covers}
    This case is trivial: If $N = L 
    \rtimes_{B} \Z$, then $L = \Z^{2}$ or $L = (A - I)\Z^{2}$.
    Since ${[\Z^{2}:(A - I)\Z^{2}]} = q$ 
    is prime, there can be no intermediate subgroups.

\subsection{Two-fold covers}  \label{sec:2-fold}
    This case is also fairly simple.  The index
    ${[\Z^{2}:(A^{2} - I)\Z^{2}]}$ is equal to $3rq$, so 
    the intermediate subgroups are of index $3$, $q$, $3q$, $rq$, $3r$, and 
    $r$.  We will treat the first three, and consider the rest as 
    dual to these.
    
    Since the index of each subgroup in this case is a prime or a product 
    of distinct primes, each subgroup is uniquely determined by its index.
    From the one-fold case, we already know that the index $q$ 
    subgroup is $(A-I)\Z^{2}$ and does not change the action 
    of $g$.  The dual subgroup is then $(A+I)\Z^{2}$ and also 
    does nothing to $g$.
    
    Consider the index 3 subgroup as the kernel of a map 
    $\Z^{2} \to \Z_{3}$.  There are four distinct such
    mappings; their kernels are
    \begin{eqnarray}
	\{(x,y) \in \Z^{2} &|& y \equiv 0\pmod 3 \}, \nonumber \\
	\{(x,y) \in \Z^{2} &|& x \equiv 0\pmod 3 \}, \nonumber \\
	\{(x,y) \in \Z^{2} &|& x+y \equiv 0\pmod 3 \},\;\mathrm{and} \nonumber \\
	\{(x,y) \in \Z^{2} &|& x-y \equiv 0\pmod 3 \}.\label{index-3}
    \end{eqnarray}
    Modulo 3, our matrix $A$ is $\left(\begin{smallmatrix} 2&2 \\ 0&2 
    \end{smallmatrix}\right)$, and it is easy to check that only the first 
    subgroup in (\ref{index-3}) is preserved by $A$.  The index 3 
    normal subgroup, then, is given by
    \begin{eqnarray*}
	\mathbf{3} &=& \langle \left(\begin{smallmatrix} 1 \\ 0 
	\end{smallmatrix}\right),\, 
	\left(\begin{smallmatrix} 0 \\ 3 \end{smallmatrix}\right) \rangle\,.
    \end{eqnarray*}
    Conjugating $A$ by $\left(\begin{smallmatrix} 1&0 \\ 0&3 
    \end{smallmatrix}\right)$ gives the action of $g$ on this 
    subgroup; it is
    \begin{eqnarray*}
	A_{\mathbf{3}} = \left(\begin{array}{cc} 1640 & 663 \\ -47 & -19 
	\end{array}\right).
    \end{eqnarray*}
    
    Knowing the subgroups of index 3 and $q$, one can find the subgroup 
    of index $3q$ by acting on the index 3 subgroup with $(A-I)$.  The 
    action of $g$ on this subgroup is then given by 
    $\left[(A-I)\left(\begin{smallmatrix} 1&0 \\ 0&3 
    \end{smallmatrix}\right)\right]^{-1} A \left[(A-I)\left(\begin{smallmatrix} 1&0 \\ 0&3 
    \end{smallmatrix}\right)\right]$, which is just $A_{\mathbf{3}}$.
    
    Calculating the cycle and dual cycle of $A_{\mathbf{3}}$, we find 
    they have lengths 4 and 42, so the monodromy $(A_{\mathbf{3}})^{2}$ has 
    cycle and dual cycle of lengths 8 and 84, well over the limit for 
    complete intersections.
    
    The index $r$ subgroup represents the cover dual to the one 
    represented by the $3q$ subgroup, but we will need a basis 
    representation for it 
    in calculating four-fold covers.  As in the index 3 case, the 
    possible subgroups are generated by $\{(1,0),\,(0,r)\}$ or 
    $\{(t,1),\,(r,0)\}$, for $0\leq t < r$.  Only the second of these is 
    preserved by $A$.  Conjugating by 
    $\left(\begin{smallmatrix} t&r \\ 1&0 \end{smallmatrix}\right)$ 
    yields
    \begin{eqnarray*}
	\left(\begin{array}{cc} -19-141t&-141r \\ 
	\frac{1}{r}(221+1659t+141t^{2})&1640+141t \end{array}\right) ,
    \end{eqnarray*}
    which suggests we solve
    \begin{eqnarray}
	221+1659t+141t^{2} \equiv 0\pmod r . \label{eqn:solve-t}
    \end{eqnarray}
    For $r = 541$, we find $t = 138$ is the unique solution 
    (modulo $r$), so the index $r$ subgroup is
    \begin{eqnarray*}
	\mathbf{r} &=& \left\langle \left(\begin{smallmatrix} 138 \\ 1 
	\end{smallmatrix}\right),\, \left(\begin{smallmatrix} 541 \\ 0 
	\end{smallmatrix}\right)\right\rangle .
    \end{eqnarray*}

    \subsection{Three-fold covers} \label{sec:3-fold}
    Here we have ${[\Z^{2}:(A^{3} - I)\Z^{2}]} = q(2s)^{2}$, so the group 
    structure is somewhat more complicated: we have to contend with subgroups 
    of index $s^{2}$, as well as subgroups of composite index.  Here 
    we will look at the subgroups of index $q$, 2, $s$, $2s$, $2q$, 
    $sq$, $s^{2}$, and $s^{2}q$, the rest being dual to these.  
    The index $q$ subgroup is still $(A-I)\Z^{2}$, and 
    we can use the same strategy employed in the two-fold case to 
    move from index $x$ to index $xq$, i.e., multiply by $(A-I)$.  
    Thus only the index $2$, $s$, and $s^{2}$ cases require explicit 
    calculation.

    There are three mappings onto $\Z_{2}$, with kernels
    \begin{eqnarray}
	\{(x,y) \in \Z^{2} &|& y \equiv 0\pmod 2 \}, \nonumber \\
	\{(x,y) \in \Z^{2} &|& x \equiv 0\pmod 2 \},\;\mathrm{and} \nonumber \\
	\{(x,y) \in \Z^{2} &|& x+y \equiv 0\pmod 2 \}.\label{index-2}
    \end{eqnarray}
    But $A \equiv \left(\begin{smallmatrix} 0&1 \\ 1&1 
    \end{smallmatrix}\right)\pmod 2$, so none of these subgroups are 
    preserved.  Therefore there is no normal cover of index 2 in the 
    fiber.
        
    Using the fact that there are no index 2 normal subgroups, a 
    simple argument eliminates the $2s$ and $2q$ cases as well:  
    Suppose one of these latter is normal; let $K$ be this normal 
    fiberwise subgroup of $G$.  But then $|G/K| = 2s$ or $2q$, and in 
    either case $G/K$ has a characteristic subgroup of order 2, so 
    there must be a normal index 2 subgroup of $G$.  We already know 
    this is not the case.
    
    There are $s+1$ possibilities for index $s$ subgroups:
    \begin{eqnarray}
	\{(x,y) \in \Z^{2} &|& y \equiv 0\pmod s \},\; \mathrm{and} \nonumber \\
	\{(x,y) \in \Z^{2} &|& x-ty \equiv 0\pmod s ,\;0\leq t <s \}.
    \end{eqnarray}
    The requirement that $A$ preserve these subgroups is equivalent 
    to asking that $P^{-1}AP \in SL_{2}(\Z)$, where $P$ is a 
    change of basis matrix.  We can easily check that 
    $\left(\begin{smallmatrix} 1&0 \\ 0&s \end{smallmatrix}\right)^{-1} 
    \left(\begin{smallmatrix} 1640&221 \\ -141&-19 \end{smallmatrix}\right) 
    \left(\begin{smallmatrix} 1&0 \\ 0&s \end{smallmatrix}\right)$ 
    does not have integer entries, so the first subgroup is ruled 
    out.  For the other $s$ possibilities, conjugation by 
    $\left(\begin{smallmatrix} t&r \\ 1&0 \end{smallmatrix}\right)$ 
    gives
    \begin{eqnarray*}
	\left(\begin{array}{cc} -19-141t&-141s \\ 
	\frac{1}{s}(221+1659t+141t^{2})&1640+141t \end{array}\right),
    \end{eqnarray*}
    so as in (\ref{eqn:solve-t}), we solve $221+1659t+141t^{2} \equiv 
    0\pmod s$.  With $s=811$, we find $t=183$ and $t=668$ are the two solutions.
    
    The two index $s$ normal subgroups are then
    \begin{eqnarray*}
	\mathbf{s_{1}} = \left\langle\left(\begin{smallmatrix} 183\\ 1 
	\end{smallmatrix}\right), 
	\left(\begin{smallmatrix} 811\\ 0 \end{smallmatrix}\right)\right\rangle ; & &
	\mathbf{s_{2}} = \left\langle\left(\begin{smallmatrix} 668\\ 1 
	\end{smallmatrix}\right), 
	\left(\begin{smallmatrix} 811\\ 0 \end{smallmatrix}\right)\right\rangle ,
    \end{eqnarray*}
    and the action of $g$ on each is
    \begin{eqnarray*}
	A_{\mathbf{s_{1}}} &=& \left(\begin{array}{cc} -25822 & -114351 \\ 
	6197 & 27443 \end{array}\right)\,; \nonumber \\
	A_{\mathbf{s_{2}}} &=& \left(\begin{array}{cc} -94207 & -114351 \\ 
	78947 & 95828 \end{array}\right)\,.
    \end{eqnarray*}
    We now calculate the relevant cycles, and find that the cycle and 
    dual cycle of $(A_{\mathbf{s_{1}}})^{3}$ have lengths $3 \cdot 149$ 
    and $3 \cdot 3$, respectively; while the cycle and dual of 
    $(A_{\mathbf{s_{2}}})^{3}$ have lengths $3 \cdot 14$ and $3 \cdot 12$.  None 
    of these, then, is a complete intersection.
    
    From the Sylow theorems we know there is exactly one order 
    $s^{2}$ subgroup of $\Z^{2}/(A^{3}-I)\Z^{2}$, so 
    there can be only one normal subgroup of $\Z^{2}$ of 
    index $s^{2}$.  We can easily check that this subgroup is
    \begin{eqnarray*}
	\mathbf{s^{2}} &=& \left\langle\left(\begin{smallmatrix} s\\ 0 
	\end{smallmatrix}\right), 
	\left(\begin{smallmatrix} 0\\ s \end{smallmatrix}\right)\right\rangle .
    \end{eqnarray*}
    Since $s$ divides $(A^{3}-I) = (A-I)(A^{2}+A+I)$, 
    $(A^{3}-I)\Z^{2} \subset \left(\begin{smallmatrix} s&0 \\ 
    0&s \end{smallmatrix}\right)\Z^{2}$.  Also, this subgroup 
    is clearly preserved by $A$, since
    \begin{eqnarray*}
	A_{\mathbf{s^{2}}} &=& \left(\begin{array}{cc} s&0 \\ 0&s 
	\end{array}\right)^{-1} A \left(\begin{array}{cc} s&0 \\ 0&s 
	\end{array}\right) = A.
    \end{eqnarray*}
    
    We are now left with the $sq$ and $s^{2}q$ cases.  These subgroups are 
    found simply by acting with $(A-I)$ on the bases of the index $s$ and index 
    $s^{2}$ subgroups.  If $P$ represents the change of basis from 
    one of these to $\Z^{2}$, then $(A-I) P$ is the 
    change of basis for the corresponding subgroup with the added 
    factor of $q$.  Then 
    $\left[(A-I)P\right]^{-1}A\left[(A-I)P\right] = 
    P^{-1}(A-I)^{-1}A(A-I)P = P^{-1}AP$, so we have
    \begin{eqnarray*}
	A_{\mathbf{(sq)_{\mathit{i}}}} &=& A_{\mathbf{s_{\mathit{i}}}}\;(i=1,2);\nonumber \\
	A_{\mathbf{(s^{2}q)}} &=& A_{\mathbf{s^{2}}}.
    \end{eqnarray*}
    We have already examined both of these matrices, so we are done 
    with three-fold covers.

\subsection{Four-fold covers} \label{sec:4-fold}
    This case is similar to the three-fold case, but simplified by 
    the fact that we have already treated most of the relevant 
    subgroups.  The index is now 
    $3p^{2}qr$, and we can use the same duality argument to focus our 
    attention on half of the subgroups; of these, only those of index 
    $p$, $qp$, $rp$, and $3p$ have not been discussed above.
    
    As in the three-fold case, the $p+1$ subgroups of index $p$ are
    \begin{eqnarray}
	\{(x,y) \in \Z^{2} &|& y \equiv 0\pmod p \}\; \mathrm{and} \nonumber \\
	\{(x,y) \in \Z^{2} &|& x-ty \equiv 0\pmod p ,\;0\leq t < p \}.
    \end{eqnarray}
    Again, the first of these can be immediately eliminated, since it 
    is not preserved by $A$.  We then solve as before for $t$ such that 
    $221+1659t+141t^{2} \equiv 0\pmod p$, finding $t=139$ or $t=608$.  
    
    Thus, the two index $p$ normal covers are
    \begin{eqnarray*}
	\mathbf{p_{1}} = \left\langle\left(\begin{smallmatrix} 139\\ 1 
	\end{smallmatrix}\right), 
	\left(\begin{smallmatrix} 1621\\ 0 \end{smallmatrix}\right)\right\rangle ; & &
	\mathbf{p_{2}} = \left\langle\left(\begin{smallmatrix} 608\\ 1 
	\end{smallmatrix}\right), 
	\left(\begin{smallmatrix} 1621\\ 0 \end{smallmatrix}\right)\right\rangle ,
    \end{eqnarray*}
    and the action of $g$ on each is
    \begin{eqnarray*}
	A_{\mathbf{p_{1}}} &=& \left(\begin{array}{cc} -19618 & -228561 \\ 
	1823 & 21239 \end{array}\right)\,; \nonumber \\
	A_{\mathbf{p_{2}}} &=& \left(\begin{array}{cc} -85747 & -228561 \\ 
	32777 & 87368 \end{array}\right)\,.
    \end{eqnarray*}
    The cycle and 
    dual cycle of $(A_{\mathbf{p_{1}}})^{4}$ have lengths $4 \cdot 12$ 
    and $4 \cdot 10$, respectively; the cycle and dual of 
    $(A_{\mathbf{p_{2}}})^{4}$ have lengths $4 \cdot 101$ and $4 
    \cdot 3$.  None of these is a CI.
    
    We can find the composite-index subgroups using a method similar to the 
    one used in the three-fold case.  As in that case, 
    the action does not change as we move from the index $p$ to the 
    index $qp$ subgroups:
    \begin{eqnarray*}
	A_{\mathbf{(qp)_{\mathit{i}}}} = A_{\mathbf{p_{\mathit{i}}}}.
    \end{eqnarray*}
    
    Using the representation of the index $r$ subgroup found in 
    Section \ref{sec:2-fold}, we can calculate the normal index $rp$ 
    subgroups.  Now we want
    \begin{eqnarray*}
	P_{t}^{-1} A P_{t} \in SL_{2}(\Z),
    \end{eqnarray*}
    where $P_{t} = \left(\begin{smallmatrix} 138&541 \\ 1&0 
    \end{smallmatrix}\right) \left(\begin{smallmatrix} t&p \\ 1&0 
    \end{smallmatrix}\right)$.  This reduces to solving
    \begin{eqnarray}
	1527 + 1571t + 1097t^{2} \equiv 0\pmod p .
    \end{eqnarray}
    We find that $t=541$ and $t=653$ solve this, so the actions 
    of $g$ on the two index $rp$ subgroups are
    \begin{eqnarray*}
	A_{(\mathbf{rp})_{1}} = P_{541}^{-1} A P_{541} 
	= \left(\begin{array}{cc} 2935465 & 8732327 \\ 
	-986243 & -2933844 \end{array}\right); \\
	A_{(\mathbf{rp})_{2}} = P_{653}^{-1} A P_{653} 
	= \left(\begin{array}{cc} 3538809 & 8732327 \\ 
	-1433459 & -3537188 \end{array}\right),
    \end{eqnarray*}
    where $P_{541} = \left(\begin{smallmatrix} 138&541 \\ 1&0 
    \end{smallmatrix}\right) \left(\begin{smallmatrix} 541&1621 \\ 1&0 
    \end{smallmatrix}\right)$
    and $P_{653} = \left(\begin{smallmatrix} 138&541 \\ 1&0 
    \end{smallmatrix}\right) \left(\begin{smallmatrix} 653&1621 \\ 1&0 
    \end{smallmatrix}\right)$.  
    The classifying matrix $(A_{\mathbf{(rp_{1})}})^{4}$ has cycle length 
    $4 \cdot 8$ and dual cycle 
    length $4 \cdot 10$, and $(A_{\mathbf{(rp)_{2}}})^{4}$ has cycle length 
    $4 \cdot 8$ and dual cycle length $4 \cdot 34$; the index $rp$ covers, 
    then, are not complete intersections.
    
    Finally, we turn to the index $3p$ subgroups.  In fact, we have 
    already calculated the monodromies and cycles for these.  Recall that adding a 
    factor of $q$ to the index leaves the monodromy unchanged; thus 
    the index $3p$ covers have the same 
    monodromies as the index $3pq$ covers.  But the index $3pq$ covers are 
    dual to the index $rp$ covers, whose classifying matrices (and 
    duals) were computed above.  This completes our verification that 
    the cusp classified by $A$, with cycle $(8,2,4,3,12)$, has no 
    Galois cover by a complete intersection.

\bibliography{galcov}
\bibliographystyle{amsplain}

\small{\textsc{Department of Mathematics, Columbia University, New York, NY 
10027}

\textit{E-mail address}: \texttt{anderson@math.columbia.edu}}

\end{document}